\documentstyle[12pt]{article}
\newlength{\pubnumber} 

\catcode`\@=11
\@addtoreset{equation}{section}

\def\section{\@startsection{section}{1}{\z@}{3.5ex plus 1ex minus .2ex}
 {2.3ex plus .2ex}{\large\bf}}
\def\subsection{\@startsection{subsection}{2}{\z@}{2.3ex plus .2ex}
 {2.3ex plus .2ex}{\bf}}

\begin{document}

\begin{titlepage}
\samepage{
\setcounter{page}{1}
\vskip 5.5truecm
\begin{center}
  {\large\bf Character States and Generator Matrix
   Elements for $Sp(4) \supset SU(2) \times U(1)$\\}

 \vskip 1.5truecm
     {\bf N. Hambli, J. Michelson and R. T. Sharp \\}
   {Department of Physics, McGill University\\
    3600 University St., Montr\'eal, Qu\'ebec~H3A-2T8~~Canada}
\footnote{Research supported in part by the Natural Sciences
and Engineering Research Council of Canada and by the Fonds
FCAR du Qu\'ebec.}
\end{center}
\vfill\eject
\begin{abstract}
   {A new set of polynomial states (to be called character
    states) are derived for $Sp\,(4)$ reduced to its
    $SU(2) \times U(1)$ subgroup, and the relevant generator
    matrix elements are evaluated for generic representations
    $(a,b)$ of $Sp(4)$. (The degenerate representations
    $(a,0)$ and $(0,b)$ were treated in our previous work
    and are also given in this paper). The group--subgroup in
    question is that of the seniority model of nuclear physics.\\}
\end{abstract}
 }
\end{titlepage}

\setcounter{footnote}{0}
\def\beq{\begin{equation}}
\def\eeq{\end{equation}}
\def\beqn{\begin{eqnarray}}
\def\barrl{\begin{array}{ll}}
\def\eeqn{\end{eqnarray}}
\def\earr{\end{array}}


\section{Introduction}

The group--subgroup $Sp(4) \supset SU(2) \times U(1)$
finds application in the nuclear seniority model
[1]. Polynomial basis states (in the states
of the $Sp(4)$ fundamental representations) have been
given by Hecht [2], Parikh [3], Ahmed and Sharp
[4] and by Smirnov and Tolstoy [5]. In this paper
we present new polynomial basis states, which we christen
``character states''; we believe them to be simpler than
any of the earlier ones. Analogous states have been given
recently for $SU(3)$ and $SO(5)$ reduced according to
their finite Demazure--Tits subgroup [6], for $SO(7)$
reduced according to $SU(2)^3$ [7] and for
$G_2$ reduced according to $SU(3)$ [8,9].

In Section~(II) we discuss character states in the
context of $Sp(4) \supset SU(2) \times U(1)$.
Section~(III) deals with generator matrix elements
for the degenerate representations $(a,0)$ and
$(0,b)$. Section~(IV) treats generic representations.
Section~(V) contains some concluding remarks.

\section{Character states for $Sp(4) \supset SU(2) \times U(1)$}

The basis states of the $Sp(4)$ representation
$(a,b)$ are polynomials of degrees $a,b$ in the
states of the respective fundamental representations.
Hence only stretched IR's (representation labels
additive) in the direct product of $a$ copies of $(1,0)$
and $b$ copies of $(0,1)$ are to be retained.

Consider the quadratic direct products
\begin{eqnarray}
{(1,0)^2}_{10} & = & (2,0)_{10}\; ,\nonumber\\
{(0,1)^2}_{15} & = & (0,2)_{14} + (0,0)_{1}\; ,
\label{aa}\\
(1,0) \times (0,1)_{20} & = & (1,1)_{16} + (1,0)_{4}\; .
\nonumber
\end{eqnarray}
The square of a representation above means the symmetric
(polynomial) part of the direct product of two copies.
A subscript on a representation or product is its dimension.
The stretched, or wanted, part of each product is the first
representation on the right; the states of the other
representations are unwanted for the purpose of forming
our polynomial basis.

The $Sp(4)$ character generator may be written
\begin{equation}
X(A,B;M,Z) =
{1\over{\beta \eta \zeta}}\;
\left[{1\over{\xi \alpha \theta}} +
{\gamma\over{\alpha \theta \gamma}} +
{\kappa\over{\theta \gamma \kappa}} +
{\delta\over{\gamma \kappa \delta}} +
{{\alpha\delta}\over{\gamma \delta \alpha}}\right]\; ;
\label{ab}
\end{equation}
we have adopted the space--saving convention that a
variable in a denominator stands for unity minus that
variable. To stress our interpretation of the character
generator as generating function for basis states we have
expressed $X$ in terms of the fundamental representation
states (see figure 1). In the power series expansion of
(\ref{ab}) each term of degree $a$ in the $(1,0)$
variables and degree $b$ in the $(0,1)$ variables
represents one state of the representation $(a,b)$.
The states thus defined are complete and non
redundant. We call them character states. We remark
that they are contaminated by unwanted states belonging
to lower representations; that does not matter for the
purpose of computing generator matrix elements. For $X$
as a generating function for characters, or weights,
the following substitutions should be made:
\begin{eqnarray}
\alpha & \to & A M^{1\over2} Z^{1\over2},\quad
\beta  \to  A M^{1\over2} Z^{-\, {1\over2}},\quad
\gamma  \to  A M^{-{1\over2}} Z^{1\over2}, \quad
\delta  \to  A M^{-{1\over2}} Z^{-{1\over2}}\; ,
\nonumber\\
\eta & \to & B Z,\quad
\xi  \to  B M~,\quad
\zeta  \to  B Z^{-1}, \quad
\kappa  \to  B M^{-1}~, \quad\
\theta  \to  B\; .
\label{ac}
\end{eqnarray}
The variables $A$ and $B$ are dummies which carry the
$Sp(4)$ representation labels $a$ and $b$ as exponents;
$M$ carries the $SU(2)$ weight $m$ and $Z$ the $U(1)$
label $z$. Then in the expansion of $X$,
\begin{eqnarray}
X(A,B;M,Z) = \sum_{abmz}\; A^a B^b M^m Z^z C_{abmz}\; ,
\label{ad}
\end{eqnarray}
the coefficient $C_{abmz}$ is the multiplicity of the
weight $(m,z)$ in the IR $(a,b)$. The character generator
(\ref{ab}) with the substitutions (\ref{ac}) agrees with
earlier versions [10,11] when the weights
are expressed in the same basis and the terms put over a
common denominator.

Examination of (\ref{ab}) reveals that certain pairs
of variables, namely $\alpha\kappa$, $\gamma\xi$,
$\delta\xi$, $\delta\theta$, $\xi\kappa$ never
appear in the same term and therefore these
products never appear in the expression for a
state; we say they are incompatible.
Each incompatible pair appears as one term in
the expression for one of the unwanted states
on the right hand side of (\ref{aa}).
Setting the unwanted states equal to zero
and solving for the incompatible pairs gives
the following substitutions by which incompatible
pairs may be eliminated in favour of compatible ones
when they arise in the course of a calculation:
\begin{eqnarray}
\alpha\, \kappa & = &  \delta \eta +
{\sqrt{2}\over 2}\; \gamma \theta\; ,
\nonumber\\
\gamma\, \xi  & = & -\; \beta \eta +
{\sqrt{2}\over 2}\; \alpha \theta\; ,
\nonumber\\
\delta\, \xi  & = &  \alpha \zeta +
{\sqrt{2}\over 2}\; \beta \theta\; ,
\label{ae}\\
\delta\, \theta & = &  \sqrt{2}\; \beta \kappa +
\sqrt{2}\; \gamma \zeta\; ,
\nonumber\\
\xi\, \kappa & = & \eta \zeta + {1\over 2}\; \theta^2\; .
\nonumber
\end{eqnarray}

It is shown by Farell, Lam and Sharp [8], and Hambli and Sharp
[9] that elementary unwanted states are all of degree $2$ and
hence that elementary incompatibilities are between pairs of
states only, a fact verified straightforwardly for $Sp(4)$.

We complete this section by giving the
$Sp(4) \supset SU(2) \times U(1)$ branching
rules generating function:
\begin{equation}
{1\over{\left(1 - A T^{1\over 2} Z^{1\over 2}\right)
\left(1 - A T^{1\over 2} Z^{- {1\over 2}}\right)
\left(1 - B Z\right)
\left(1 - B Z^{-1}\right)}}
\left[
{1\over{\left(1 -  B T\right)}} +
{{A^2}\over{\left(1 - A^2\right)}}
\right]\; .
\label{af}
\end{equation}
The elements of (\ref{af}) may be interpreted as
the highest states of subgroups IR's contained in
low $Sp(4)$ IR's:
\begin{eqnarray}
\alpha & \sim & A T^{1\over2} Z^{1\over2},\quad
\beta  \sim  A T^{1\over2} Z^{-\, {1\over2}},\quad
\alpha\delta - \beta\gamma  \sim  A^2 \; ,
\nonumber\\
\eta & \sim & B Z,\quad
\zeta  \sim  B Z^{- 1},\quad
\xi  \sim  B T \; .
\label{ag}
\end{eqnarray}
Thus the highest state of a subgroup IR
may be written, in unnormalized form, and
representation labels suppressed,
\begin{equation}
\left|
\, t\;,t\;,z\; ; v\,
\right\rangle =
\alpha^x \; \beta^y\; \eta^u\;
\zeta^v\; \xi^w\;
\left(\alpha\delta - \beta\gamma\right)^s\; ,
\label{ah}
\end{equation}
with
\beqn
a & = & x + y + 2s ~ , \quad
b  =  u + v + w \; ,
\nonumber\\
& &\label{ai}\\
t & = & {1\over 2}\; \left(x + y\right) + w~ , \quad
z   =  {1\over 2}\; \left(x - y\right) + u - v \; ,
\nonumber
\eeqn
and $s\; w = 0$. We allow $v$ to play the role of the
missing label.

We may distinguish two types of state, according
to whether (type I) $s = 0$ (and $t \ge a/2$)
or (type II) $w = 0$ (and $t \le a/2$). Solving
(\ref{ai}) for $x$, $y$, $u$ and $w$ or $s$,
(\ref{ah}) becomes
\begin{equation}
\left|
\,t\;,t\;,z\; ; v\,
\right\rangle =
\alpha^{t + z - b + 2 v} \; \beta^{a + b - t - z - 2 v} \;
\eta^{{a\over 2} + b - t - v} \;
\zeta^v\;
\xi^{t - {a\over 2}} \; ,
\label{aj}
\end{equation}
for type I, and, for type II
\begin{equation}
\left|
\,t\;,t\;,z\; ; v\,
\right\rangle =
\alpha^{t + z - b + 2 v} \;
\beta^{b + t - z - 2 v}\;
\left(\alpha\delta - \beta\gamma\right)^{{a\over 2} - t}\;
\eta^{b - v}\;
\zeta^v\; .
\label{ak}
\end{equation}
For $t = a/2$ the two types of state coincide.

The $SU(2)$ root generators are
\begin{eqnarray}
T_{+} & = &  \alpha\, \partial_{\gamma} +
\beta \partial_{\delta} +
\sqrt{2}\, \left(\xi\, \partial_\theta +
\theta\, \partial_\kappa \right)\; ,
\nonumber\\
& & \label{al}\\
T_{-}  & = &  \gamma\, \partial_{\alpha} +
\delta \partial_{\beta} +
\sqrt{2}\, \left(\theta\, \partial_\xi +
\kappa\, \partial_\theta \right)\; .
\nonumber
\end{eqnarray}
The other $6$ $Sp(4)$ root generators
constitute an $SU(2)$ vector $G$ with
$z = 1$,
\begin{eqnarray}
G_{+1} & = &  \alpha\, \partial_\delta +
\eta\, \partial_\kappa + \xi\, \partial_\zeta \; ,
\nonumber\\
G_{-1}  & = &  -\; \gamma\, \partial_\beta +
\kappa \partial_\zeta +
\eta\, \partial_\xi \; ,
\label{am}\\
G_{0} & = & {\sqrt{2}\over 2}\,
\left(\gamma\, \partial_\delta -
\alpha\, \partial_\beta \right) +
\theta\, \partial_\zeta -
\eta\, \partial_\theta \; ,
\nonumber
\end{eqnarray}
and a vector $\overline{G}$ with $z = - 1$,
\begin{eqnarray}
\overline{G_{+1}} & = &  \beta\, \partial_\gamma -
\xi\, \partial_\eta - \zeta\, \partial_\kappa \; ,
\nonumber\\
\overline{G_{-1}}  & = &  -\; \delta\, \partial_\alpha -
\kappa \partial_\eta - \zeta\, \partial_\xi \; ,
\label{an}\\
\overline{G_{0}} & = & {\sqrt{2}\over 2}\,
\left(\delta\, \partial_\gamma  -
\beta\, \partial_\alpha \right) +
\zeta\, \partial_\theta -
\theta\, \partial_\eta \; .
\nonumber
\end{eqnarray}
$G$ and $\overline{G}$ are hermitian conjugate  with
\begin{eqnarray}
\overline{G_{i}} = \left(- 1\right)^i\;\;
G_{- i}^{\dagger}\; .
\label{ao}
\end{eqnarray}

We will save considerable work later by noticing
that under the substitutions
$\alpha \leftrightarrow \delta$,
$\beta \leftrightarrow \gamma$,
$\eta \leftrightarrow \zeta$,
$\xi \leftrightarrow \kappa$,
we have
$G_i \leftrightarrow -\, \overline{G_{-i}}$
and that under the same substitutions we have
for the type I states (\ref{aj})
$\left|\,t\;, m\;, z\; ; v\,\right\rangle
\leftrightarrow
\left|\,t\;, -m\;, -z\; ; a/2 + b - t - v\,\right\rangle$
while for type II states (\ref{ak})
$\left|\,t\;, m\;, z\; ; v\,\right\rangle
\leftrightarrow
\left|\,t\;, -m\;, -z\; ; b - v\,\right\rangle$.

\section{Generator matrix elements for degenerate
representations}

For the degenerate IR's $(a,0)$ and $(0,b)$ there
is no missing label and the basis states are orthogonal
and hence can be normalized straightforwardly. Since
in that sense they can be handled more satisfactorily
than generic states $(a > 0 , b > 0)$ we treat them
separately.

For $(a,0)$, according to (\ref{ak}) with $v = 0$, the
highest state of an $SU(2)$ IR is (we suppress the
label $a$)
\begin{equation}
\left|
\,t\;,t\;,z\; \,
\right\rangle = N_{t\;\;z}\;\;
\alpha^{t + z}\; \beta^{t - z}\;
\left(\alpha\delta - \beta\gamma\right)^{{a\over 2} - t}\; .
\label{ba}
\end{equation}
The branching rule is $a/2 \ge t \ge |z|$ with
$2 t$ and $2 z$ having the parity of $a$.

$(a,0)$ states are not contaminated by unwanted states
and can be normalized by standard methods.
For brevity write
\begin{equation}
\left|\, h \,\right\rangle = N_{h}\;\;
\alpha^{f}\; \beta^{g}\;
\left(\alpha\delta - \beta\gamma\right)^{h}\; ,
\label{bb}
\end{equation}
and equate
\begin{equation}
\left\langle\, h + 1\,\right|\,
\alpha\, \delta - \beta \,\gamma \,
\left|\, h \,\right\rangle =
\left\langle\, h \,\right|\,
\partial_\alpha \, \partial_\delta -
\partial_\beta \, \partial_\gamma \,
\left|\, h + 1\,\right\rangle \; ,
\label{bc}
\end{equation}
to obtain a recursion relation for $N_h$,
whose solution is, for $N_{t\;\;z}$,
\begin{equation}
N_{t\;\;z} =
\sqrt{
{{(2 t + 1)!}
\over
{\left(t + z\right) !\; \left(t - z\right) !\;
\left({a/2} - t\right) !\;
\left({a / 2} + t + 1 \right) !}}
}\; .
\label{bd}
\end{equation}

We can write
\begin{equation}
G_0 \; \left|\, t\; , t\; , z\,\right\rangle = \;
A\; \left|\, t+1\; , t\; , z+1\,\right\rangle +
B\; \left|\, t\; , t\; , z+1\,\right\rangle \; ,
\label{bd}
\end{equation}
where the matrix elements $A$ and $B$ will now
be determined. Apply $T_{+}$ to
(\ref{bd}) with the result
\begin{equation}
\sqrt{2}\; G_1 \; \left|\, t\; , t\; , z\,\right\rangle = \;
\sqrt{2(t+1)}\; A\; \left|\, t+1\; , t+1\; , z+1\,\right\rangle \; ,
\label{be}
\end{equation}
from which
\begin{eqnarray}
A & = &
\left\langle\, t+1\;,t\;,z+1\,\right|\,
G_0 \,
\left|\, t\;,t\;,z \,\right\rangle
= {{{a/2} - t}\over{\sqrt{t+1}}}\;\;
{{N_{t\;\;z}}\over{N_{t+1\;\;z+1}}}\; ,
\nonumber\\
& = &
{1\over{t + 1}}\;\;
\sqrt{
{{\left(t + z + 1\right) \left(t + z + 2\right)
\left({a / 2} - t\right) \left({a / 2} + t + 2\right)}
\over
{2\, \left(2 t + 3\right)}}
}\; .
\label{bf}
\end{eqnarray}

Inserting this value for $A$ in (\ref{bd}) we obtain
for $B$
\begin{eqnarray}
B & = &
\left\langle\, t\;,t\;,z+1\,\right|\,
G_0 \,
\left|\, t\;,t\;,z \,\right\rangle
= -\; {{({a / 2} + 1)(t - z)}\over{\sqrt{2} (t + 1)}}\;\;
{{N_{t\;\;z}}\over{N_{t\;\;z+1}}}\; ,
\nonumber\\
& = &
-\; {{{a/2 2} + 1}\over{\sqrt{2} (t + 1)}}\;\;
\sqrt{
\left(t + z + 1\right) \left(t - z\right)
}\; .
\label{bg}
\end{eqnarray}

To obtain a matrix element of $G_0$ in which $G_0$
reduces $t$ by unity make the following transformations:
\begin{eqnarray}
\left\langle t+1\;,t\;,z+1\right|
G_0
\left| t\;,t\;,z \right\rangle
& = & -
\left\langle t+1\;,-t\;,-z-1\right|\,
\overline{G_0}
\left|\, t\;,-t\;,-z \,\right\rangle\; ,
\nonumber\\
& = &
-
\left\langle t\;,-t\;,-z\right|
G_0
\left| t+1\;,-t\;,-z-1 \right\rangle \, .
\label{bh}
\end{eqnarray}
The first step follows from the substitutions in
the paragraph following (\ref{ao}) and the second
is a consequence of (\ref{ao}).

For the reduced matrix element of $G$ we get
\begin{eqnarray}
\left\langle \,t+1\;,z+1\,
\left|\!\left| \,G\, \right|\!\right|
\,t\;,z\, \right\rangle & = &
\sqrt{
{{\left(t + z + 1\right) \left(t + z + 2\right)
\left({a / 2} - t\right) \left({a / 2} + t + 2\right)}
\over
{2 \left(t + 1\right)}}
}\; ,
\nonumber\\
\left\langle \,t\;,z+1\,
\left|\!\left| \,G\, \right|\!\right|
\,t\;,z\, \right\rangle & = &
-\; \left({a / 2} + 1\right)\;
\sqrt{
{{\left(2 t + 1\right) \left(t + z + 1\right) \left(t - z\right)}
\over
{2 t \left(t + 1\right)}}
}\;\; ,
\label{bi}\\
\left\langle \,t-1\;,z+1\,
\left|\!\left| \,G\, \right|\!\right|
\,t\;,z\, \right\rangle & = &
\sqrt{
{{\left(t - z - 1\right) \left(t - z\right)
\left({a / 2} - t + 1\right) \left({a / 2} + t + 1\right)}
\over
{2 t}}
}\;\; .
\nonumber
\end{eqnarray}
Because of (\ref{ao}) the reduced matrix elements of
$\overline{G}$ are given in terms of those of $G$
by the relation
\begin{equation}
\left\langle \,t+k\;,z-1\,
\left|\!\left| \,\overline{G}\,
\right|\!\right|\,t\;,z\, \right\rangle =
-\left(-1\right)^k \,
\left\langle \,t\;,z\,
\left|\!\left| \,G\, \right|\!\right|
\,t+k\;,z-1\, \right\rangle\; .
\label{bip}
\end{equation}

We turn to $(0,b)$ states. According to (\ref{aj})
with $v =(b - t - z)/2$ the highest state of a subgroup
multiplet is
\begin{equation}
\left| \, t\; , t \; , z \,\right\rangle =
N_{t\;\;z} \,
\eta^{{(b-t+z)}\over 2} \,
\zeta^{{(b-t-z)}\over 2} \,
\xi^t\;\; .
\label{bj}
\end{equation}
The branching rule is $b - |z| \ge t \ge 0$.
The $SU(2)$ ladder generators are now
\begin{eqnarray}
T_{+} & = &
\sqrt{2}\; \left(
\xi\, \partial_\theta + \theta\, \partial_\kappa
\right)\; ,
\nonumber\\
& &\label{bk}\\
T_{-} & = &
\sqrt{2}\; \left(
\theta\, \partial_\xi + \kappa\, \partial_\theta
\right)\; .
\nonumber
\end{eqnarray}
The components of $G$ and $\overline{G}$ are
\begin{eqnarray}
G_{+1} & = &  \eta\, \partial_\kappa +
\xi\, \partial_\zeta \; ,
\nonumber\\
G_{-1}  & = &   \kappa\, \partial_\zeta +
\eta \partial_\xi \; ,
\label{bl}\\
G_{0} & = &
\theta\, \partial_\zeta -
\eta\, \partial_\theta \; ,
\nonumber
\end{eqnarray}
\begin{eqnarray}
\overline{G_{+1}} & = &  -\, \xi\, \partial_\eta -
\zeta\, \partial_\kappa \; ,
\nonumber\\
\overline{G_{-1}}  & = &  -\; \kappa\, \partial_\eta -
\zeta \partial_\xi \; ,
\label{bm}\\
\overline{G_{0}} & = &
\zeta\, \partial_\theta -
\theta\, \partial_\eta \; .
\nonumber
\end{eqnarray}
Again $\overline{G_i} = \left(-1\right)^i \;
{G_{-i}}^\dagger$. The basis states are contaminated
by states containing as a factor the unwanted scalar
$\eta\;\zeta - \xi\;\kappa + \theta^2 /2$.

Applying  $\overline{G_{-1}}$ to
$\left|\, t\; ,t\; , z\, \right\rangle$ yields
\begin{equation}
\overline{G_{-1}} \; \left|\, t\; , t\; , z\,\right\rangle = \;
A\; \left|\, t+1\; , t-1\; , z-1\,\right\rangle +
C\; \left|\, t-1\; , t-1\; , z-1\,\right\rangle \; .
\label{bn}
\end{equation}
($t$ changes only by $\pm 1$ because it has the
parity of $b+z$.) We must now find the matrix elements
$A$ and $C$. Application of $T_{+}^2$ to
(\ref{bn}) gives
\begin{equation}
\overline{G_{1}} \; \left|\, t\; , t\; , z\,\right\rangle = \;
A\; \sqrt{(t+1)(2t+1)} \;
\left|\, t+1\; , t+1\; , z-1\,\right\rangle \; ,
\label{bo}
\end{equation}
from which
\begin{eqnarray}
A & = &
\left\langle\, t+1\;,t-1\;,z-1\,\right|\,
\overline{G_{-1}} \,
\left|\, t\;,t\;,z \,\right\rangle \; ,
\nonumber\\
& = & -\, {1\over 2}\;
{{\left(b-t+z\right)}\over{\sqrt{(t+1)(2t+1)}}}\;\;
{{N_{t\;\;z}}\over{N_{t+1\;\;z-1}}}\; .
\label{bp}
\end{eqnarray}
Then (\ref{bn}) yields immediately
\begin{eqnarray}
C & = &
\left\langle\, t-1\;,t-1\;,z-1\,\right|\,
\overline{G_{-1}} \,
\left|\, t\;,t\;,z \,\right\rangle \; ,
\nonumber\\
 & = & -\,
{{t(b+t+z+1)}\over{2t+1}}\;\;
{{N_{t\;\;z}}\over{N_{t-1\;\;z-1}}}\; .
\label{bq}
\end{eqnarray}
Applying $G_{1}$ to $\left|\,t-1\;,t-1\;,z-1\,\right\rangle$
gives
\begin{equation}
\left\langle\, t\;,t\;,z\,\right|\,
{G_{1}} \,
\left|\, t-1\;,t-1\;,z-1 \,\right\rangle
= {1\over 2}\;(b - t - z + 2)\;\;
{{N_{t-1\;\;z-1}}\over{N_{t\;\;z}}}\;\; .
\label{br}
\end{equation}
Because $\overline{G_{-1}} = -\; {G_1}^\dagger$ we can
equate the right hand side of (\ref{br}) to the negative
of the right hand side of (\ref{bq}), getting
a recursion relation for the normalization constant
\begin{equation}
{{N_{t\;\;z}}\over{N_{t-1\;\; z-1}}}  =
\sqrt{
{{(b-t-z+2)(2t+1)}\over{2(b+t+z+1)t}}
}\; \; ,
\label{bs}
\end{equation}
whose solution is
\begin{equation}
N_{t\;\;z} =
\phi \;\;
\sqrt{
{(2t+1)!!}
\over
{2^{{1\over 2}\, (b+t-z)} \;
t!\;\;
\left({{b-t-z}\over2}\right)!\;\;
(b+t+z+1)!!}
}\;\; ,
\label{bt}
\end{equation}
where $\phi$ is constant when $t$ and $z$
are given the same increment. Increasing $t$
and $z$ each by $(b-t-z)/2$ takes us to the state
\begin{equation}
\left|
{1\over 2} (b+t-z)\;,
{1\over 2} (b+t-z)\;,
{1\over 2} (b-t+z)
\right\rangle =
N_{{1\over 2} (b+t-z)\;\; {1\over 2}\, (b-t+z)}\;
\eta^{{1\over 2} (b-t+z)}\;
\xi^{{1\over 2} (b+t-z)}\; ,
\label{bu}
\end{equation}
which is on the boundary of the $(0,b)$ weight
diagrams and uncontaminated. It is easily normalized:
\begin{equation}
N_{{1\over 2}\, (b+t-z)\;\; {1\over 2}\, (b-t+z)} =
{1\over{\sqrt{
\left({{b-t+z}\over 2}\right)! \;\;
\left({{b+t-z}\over 2}\right)!
}}}\;\; .
\label{bv}
\end{equation}
Using (\ref{bv}) in (\ref{bt}) determines
$\phi$ to be
\begin{equation}
\phi =
\sqrt{
{{2^{{1\over 2}\, (b+t-z)} \;\; (2b+1)!!}
\over
{\left({{b-t+z}\over 2}\right)! \;\;
\left(b+t-z+1\right)!!}}
}\; \;\; .
\label{bw}
\end{equation}
(\ref{bt}) then gives for $N_{t\;\;z}$
\begin{equation}
N_{t\;\;z} =
\sqrt{
{{(2t+1)!!\;\; (2b+1)!!}
\over
{t!\;\;
\left({{b-t-z}\over 2}\right)!\;\;
\left({{b-t+z}\over 2}\right)!\;\;
(b+t+z+1)!!\;\; (b+t-z+1)!!}}
}\;\;\; .
\label{bx}
\end{equation}
We have normalized the wanted part of the state
(\ref{bj}) without ever isolating it.
As a check consider (\ref{bj}) with
$b=4,t=0,z=0$, {\it i.e.,}
$N_{00}\; \eta^2 \; \zeta^2$.
The wanted part of $\eta^2\;\zeta^2$,
{\it i.e.,} the part orthogonal to
$\eta\zeta
\left(\eta \zeta - \xi \kappa +
{1/ 2}\, \theta^2 \right)$ and to
$\left(\eta \zeta - \xi \kappa +
{1/ 2}\, \theta^2 \right)^2$ is
$1/63\; \left(
15\, \eta^2 \, \zeta^2 +
8\, \xi^2 \, \kappa^2 +
2\, \theta^4 +
40\, \eta\, \zeta\, \xi\, \kappa -
20\, \eta\, \zeta\, \theta^2 -
8\, \xi\, \kappa\, \theta^2\right)$
whose norm is $20/21$ which checks with
(\ref{bx}) by which $N_{00} = \sqrt{21/20}$.

Inserting the value of the normalization constant
in (\ref{bp}) and (\ref{bq}) gives the matrix elements
explicitly. For the reduced matrix elements of $G$
we find
\begin{eqnarray}
\left\langle \,t+1\;,z-1\,
\left|\!\left| \,
{G}\, \right|\!\right|
\,t\;,z\, \right\rangle & = &
-\,\sqrt{
{1\over 2} (t + 1)\,
\left({{b - t - z}}\right)\,
(b + t + z + 3)}\; ,
\nonumber\\
&  &
\label{by}\\
\left\langle \,t-1\;,z-1\,
\left|\!\left| \,
{G}\, \right|\!\right|
\,t\;,z\, \right\rangle & = &
-\,\sqrt{
{1\over 2} t\,
\left({{b - t + z + 2}\over 2}\right)\,
(b + t - z + 1)}\;\; ,
\nonumber
\end{eqnarray}
and Eq.~(\ref{bip}) again gives the reduced matrix
elements of $\overline{G}$ in terms of those of $G$.

\section{Generator matrix elements for generic
representations}

For generic representations $(a,b)$, with $a > 0$ and
$b > 0$ it is convenient to use the states in
non-orthonormal form as given by (\ref{aj}) and
(\ref{ak}) for type I and type II states respectively.
The matrix element of a generator $X$ between two states
$\left|\, \alpha\, \right\rangle$ and
$\left|\, \beta\, \right\rangle$ will be written
$\left(\, \beta\, \right|\, X\, \left|\, \alpha\, \right)$
to imply the coefficient of $\left|\, \beta\, \right\rangle$
in $X\, \left|\, \alpha\, \right\rangle$. It is not
now the same as the overlap
$\left\langle\, \beta\, \right|\, X\,
\left|\, \alpha\, \right\rangle$ as in the degenerate
representations of Section~(II) where the states were
orthonormal. The Wigner--Eckart theorem applies equally
to this definition of matrix element and the
eigenvalues and eigenstates of a hermitian operator $H$
(suppose the Hamiltonian of a physical system is given in
terms of generators) are found by standard routines.

We can write
\begin{eqnarray}
G_{-1}\;
\left|t\;, t\;, z\;, v\,\right\rangle & = &
\sum_{v'}\; \{
\left|t+1\;, t-1\;, z+1\;, v'\,\right\rangle
\nonumber\\
& & \times
\left(\, t+1\;, t-1\;, z+1\;, v'\,\right|\,
{G_{-1}} \,
\left|\, t\;, t\;, z\;, v \,\right)
\nonumber\\
& &
\left|t\;, t-1\;, z+1\;, v'\,\right\rangle
\left(\, t\;, t-1\;, z+1\;, v'\,\right|\,
{G_{-1}} \,
\left|\, t\;, t\;, z\;, v \,\right)
\nonumber\\
& &
\left|t-1\;, t-1\;, z+1\;, v'\,\right\rangle
\nonumber\\
& & \times
\left(\, t-1\;, t-1\;, z+1\;, v'\,\right|\,
{G_{-1}} \,
\left|\, t\;, t\;, z\;, v \,\right)\}\; .
\label{ca}
\end{eqnarray}
Applying $T_{+}^2$ to both sides of (\ref{ca})
yields
\begin{eqnarray}
G_{1}\; \left|\, t\;, t\;, z\;, v \,\right\rangle & = &
\sum_{v'}\;\{
\left|t+1\;, t+1\;, z+1\;, v'\,\right\rangle
\nonumber\\
& & \times
\left( t+1\;, t-1\;, z+1\;, v'\right|
{G_{-1}}
\left| t\;, t\;, z\;, v \right)\,
\sqrt{t (2 t + 1)}\} .
\nonumber\\
& &
\label{cb}
\end{eqnarray}

The only nonzero matrix elements with $t + 1$ on the
left side are
\begin{eqnarray}
\left(\, t+1\;, t-1\;, z+1\;, v-1\,\right|\,
{G_{-1}} \,
\left|\, t\;, t\;, z\;, v \,\right) & = &
{v\over\sqrt{(1 + 2 t)\; (1 + t)}}
\nonumber
\end{eqnarray}
for type I states, and
\begin{eqnarray}
\left( t+1, t-1, z+1, v-1\right|
{G_{-1}}
\left| t, t, z, v \right) & = &
{v\over\sqrt{(1 + 2 t) (1 + t)}}
\nonumber\\
\left( t+1, t-1, z+1, v\right|
{G_{-1}}
\left| t, t, z, v \right) & = &
{{({a\over 2} - t + v)}\over\sqrt{(1 + 2 t) (1 + t)}}
\nonumber
\end{eqnarray}
for type II states.

Now subtract the $t + 1$ parts from the right side of
(\ref{ca}) and apply $T_{+}$ to what is left. The right
hand side becomes
\begin{equation}
\sum_{v'}\;
\left|t\;, t\;, z+1\;, v'\,\right\rangle
\left(\, t\;, t-1\;, z+1\;, v'\,\right|\,
{G_{-1}} \,
\left|\, t\;, t\;, z\;, v \,\right)\,
\sqrt{2 t}\; ,
\nonumber
\end{equation}
from which we find that the only nonzero matrix
elements with $t$ on the left side are
\begin{eqnarray}
\left( t\;, t-1\;, z+1\;, v-1\right|
{G_{-1}}
\left| t\;, t\;, z\;, v \right) & = &
{{v (-b + t + 2 v + z)}\over{\sqrt{2 t} (1 + t)}}
\nonumber\\
\left( t\;, t-1\;, z+1\;, v\right|
{G_{-1}}
\left| t\;, t\;, z\;, v \right) & = &
{{(1 + t + v)}\over{\sqrt{2 t} (1 + t)}}
\nonumber\\
& & \times
(-a - b + t + 2 v + z)
\nonumber
\end{eqnarray}
for type I states and
\begin{eqnarray}
\left(t\;, t-1\;, z+1\;, v-1\right|
{G_{-1}}
\left| t\;, t\;, z\;, v \right) & = &
{{v (-b + t + 2 v + z)}\over{\sqrt{2 t} (1 + t)}}
\nonumber\\
\left( t\;, t-1\;, z+1\;, v\right|
{G_{-1}}
\left| t\;, t\;, z\;, v \right) & = &
{{(1 + {a\over 2} + v)}\over{\sqrt{2 t} (1 + t)}}
\nonumber\\
& & \times
(-b - t + 2 v + z)
\nonumber
\end{eqnarray}
for type II states. Now subtract the $t$ parts
from the right side of (\ref{ca}). Only the $t - 1$
states remain and we can read that the only nonzero matrix
elements are
\begin{eqnarray}
\left( t-1\;, t-1\;, z+1\;, v-1\right|
{G_{-1}}
\left| t\;, t\;, z\;, v \right) & = &
{{v (-b + t + 2 v + z)}\over{2 t (1 + 2 t)}}\;
\nonumber\\
 & & \times
(-1 - b + t + 2 v + z)
\nonumber\\
\left( t-1\;, t-1\;, z+1\;, v\right|
{G_{-1}}
\left| t\;, t\;, z\;, v \right) & = &
{1\over{2 t (1 + 2 t)}}
(-a - b + a b + b^2 + t
\nonumber\\
& &
- 2 a t - 2 b t + 2 a b t + 2 b^2 t + 3 t^2 + 2 t^3
\nonumber\\
& &
+ 2 v - 3 a v - 4 b v + 2 a b v + 2 b^2 v + 4 t v
\nonumber\\
& &
- 6 a t v - 12 b t v + 6 t^2 v + 4 v^2 - 4 a v^2
\nonumber\\
& &
- 8 b v^2 + 16 t v^2 + 8 v^3 + z - a z - 2 b z
\nonumber\\
& &
+ 2 t z - 2 a t z - 4 b t z + 4 v z - 2 a v z
\nonumber\\
& &
- 4 b v z + 12 t v z + 8 v^2 z + z^2
\nonumber\\
& &
+ 2 t z^2 + 2 v z^2)
\nonumber\\
\left( t-1\;, t-1\;, z+1\;, v+1 \right|
{G_{-1}}
\left| t\;, t\;, z\;, v \right) & = &
{{(1 + 2 t + v)}\over{2 t (1 + 2 t)}}
\nonumber\\
& & \times
(-a - b + t + 2 v + z)
\nonumber\\
& & \times
(1 - a - b + t + 2 v + z)
\nonumber
\end{eqnarray}
for type I states and
\begin{eqnarray}
\left( t-1\;, t-1\;, z+1\;, v-1\right|
{G_{-1}}
\left| t\;, t\;, z\;, v \right) & = &
{{v (b - t - 2 v - z)}\over{2 t (1 + 2 t)}}
\nonumber\\
& & \times
(-1 - b + t + 2 v + z)
\nonumber\\
\left( t-1\;, t-1\;, z+1\;, v\right|
{G_{-1}}
\left| t\;, t\;, z\;, v \right) & = &
{{(1 + {a\over 2} + t + v)}\over{2 t (1 + 2 t)}}
\nonumber\\
& & \times
(-1 + b + t - 2 v - z)
\nonumber\\
& & \times
(-b - t + 2 v + z)
\nonumber
\end{eqnarray}
for type II states.

Notice that because $t$ does not change by more than
$\pm 1$ we do not require matrix elements of $G_{-1}$
between a type I or II state and a state of the
opposite type. A type I or II state goes to a state of
the same type or to an ambiguous state with $t = a/2$
which can be regarded as being of the same type. An
ambiguous state can go to type I or II.

We now give the reduced matrix elements, obtained
from the above ordinary matrix elements using the
Wigner--Eckart theorem
\begin{eqnarray}
\left( t+1\;, z+1\;, v-1
\right|\!\left| {G} \right|\!\left|
t\;, z\;, v \right) & = &
\sqrt{3 + 2 t} v
\nonumber\\
\left( t\;, z+1\;, v-1
\right|\!\left| {G} \right|\!\left|
t\;, z\;, v \right) & = &
\sqrt{{1 + 2 t}\over{2 t (1 + t)}} v (-b + t + 2 v + z)
\nonumber\\
\left( t\;, z+1\;, v
\right|\!\left| {G} \right|\!\left|
 t\;, z\;, v \right) & = &
\sqrt{{1 + 2 t}\over{2 t (1 + t)}}
(-a - b + t + 2 v + z)
\nonumber\\
& & \times
(1 + t + v)
\nonumber\\
\left( t-1\;, z+1\;, v-1
\right|\!\left| {G} \right|\!\left|
 t\;, z\;, v \right) & = &
{{v \left(-1 - b + t + 2 v + z\right)}\over{2 t \sqrt{1 + 2t}}}\;
\nonumber\\
& & \times
\left(-b + t + 2 v + z\right)
\nonumber\\
\left( t-1\;, z+1\;, v
\right|\!\left| {G} \right|\!\left|
 t\;, z\;, v \right) & = &
{1\over{2 t \sqrt{1 + 2 t}}}
(-a - b + a b + b^2 + t
\nonumber\\
& &
- 2 a t - 2 b t + 2 a b t + 2 b^2 t + 3 t^2 + 2 t^3
\nonumber\\
& &
+ 2 v - 3 a v - 4 b v + 2 a b v + 2 b^2 v + 4 t v
\nonumber\\
& &
- 6 a t v - 12 b t v + 6 t^2 v + 4 v^2 - 4 a v^2
\nonumber\\
& &
- 8 b v^2 + 16 t v^2 + 8 v^3 + z - a z - 2 b z
\nonumber\\
& &
+ 2 t z - 2 a t z - 4 b t z + 4 v z - 2 a v z
\nonumber\\
& &
- 4 b v z + 12 t v z + 8 v^2 z + z^2
\nonumber\\
& &
+ 2 t z^2 + 2 v z^2)
\nonumber\\
\left( t-1\;, z+1\;, v+1
\right|\!\left| {G} \right|\!\left|
 t\;, z\;, v \right) & = &
{{(1 + 2 t + v)}\over{2 t \sqrt{1 + 2 t}}}
\nonumber\\
& & \times
(-a - b + t + 2 v + z)
\nonumber\\
& & \times
(1 - a - b + t + 2 v + z)
\nonumber
\end{eqnarray}
for type I states and
\begin{eqnarray}
\left(\, t+1\;, z+1\;, v-1\,
\right|\!\left| {G} \right|\!\left|\,
 t\;, z\;, v \,\right) & = &
\sqrt{3 + 2 t}\; v
\nonumber\\
\left(\, t+1\;, z+1\;, v\,
\right|\!\left| {G} \right|\!\left|
\, t\;, z\;, v \,\right) & = &
\sqrt{3 + 2 t}\; \left({a\over 2} - t + v\right)
\nonumber\\
\left(\, t\;, z+1\;, v-1\,
\right|\!\left| {G} \right|\!\left|\, t\;, z\;, v \,\right) & = &
\sqrt{{1 + 2 t}\over{2 t (t + 1)}}\;
v\; \left(- b + t + 2 v + z\right)
\nonumber\\
\left(\, t\;, z+1\;, v\,
\right|\!\left| {G} \right|\!\left|\, t\;, z\;, v \,\right) & = &
\sqrt{{1 + 2 t}\over{2 t (1 + t)}}\;
\left(1 + {a\over 2} + v\right)\;
\nonumber\\
& & \times
(- b - t + 2 v + z)
\nonumber\\
\left(\, t-1\;, z+1\;, v-1\,
\right|\!\left| {G} \right|\!\left|\, t\;, z\;, v \,\right) & = &
{{v\; (b - t - 2 v - z)}\over{2 t \sqrt{1 + 2 t}}}\;
\nonumber\\
& & \times
(-1 - b + t + 2 v + z)
\nonumber\\
\left(\, t-1\;, z+1\;, v\,
\right|\!\left| {G} \right|\!\left|\, t\;, z\;, v \,\right) & = &
{{(1 + {a\over 2} + t + v)}\over{2 t \sqrt{1 + 2 t}}}
\nonumber\\
& & \times
(-1 + b + t - 2 v - z)
\nonumber\\
& & \times
(-b - t + 2 v + z)
\nonumber
\end{eqnarray}
for type II states.

Finally, we remark that the reduced matrix
elements of $\overline{G}$ are given in terms
of those of $G$ by the generic equivalent of
(\ref{bip})
\begin{equation}
\left( \,t+k\;,z-1\;, v'\,
\right|\!| \,\overline{G}\,
|\!\left|\,t\;,z\;, v\, \right) =
-\left(-1\right)^k \,
\left( \,t\;,z\;, v\,
\right|\!\left| \,G\, \right|\!\left|
\,t+k\;,z-1\;, v'\, \right)\; .
\label{cc}
\end{equation}

Our results can be compared with those of Section~(III)
by setting one of the $Sp(4)$ representation $a$ or $b$
labels equal to zero, omitting the `missing' label $v$,
and replacing the normalization constants $N_{t\;\; z}$
of Section~(III) by unity.

\section{Concluding remarks}

Our $B_2 \supset A_1 \times U(1)$ states are of two
types, labeled I and II respectively. There are also
``ambiguous'' states which can be regarded as being of
either type. We find that a generator never connects two
states which are definitely labeled differently, {\it i.e.,\/}
I and II. Rather a generator links two states of the same
type or a state of type I or II to an ambiguous state
(it can also link two ambiguous states). Thus, we need
only two formulas for generator matrix elements, one
for each type of state. This raises the following
consideration: Is there a group--subgroup such that
a generator links two states which are unambiguously
of different types? If so a special formula would have to
be given for the relevant generator matrix elements. We
do not know of any such group--subgroup yet.

In a future paper we hope to compute generator matrix
elements for $B_3$ in a $G_2$ basis.

\pagebreak

\section{Figure Captions}

Figure 1. States of the fundamental representations of $Sp(4)$.

\pagebreak

\section{References}

\noindent
[1].~B. H. Flowers,
{\it Proc. Roy. Soc.} (London) {\bf A212}
248 -- 263
\medskip

\noindent
[2].~K. T. Hecht,
{\it Nucl. Phys.} {\bf 63} 177 -- 213
\medskip

\noindent
[3].~J. C. Parikh, {\it Nucl. Phys.} {\bf 63}
214 -- 232
\medskip

\noindent
[4].~K. Ahmed and R. T. Sharp,
{\it J. Math. Phys.} {\bf 11}
1112 -- 1117
\medskip

\noindent
[5].~Y. F. Smirnov and V. N. Tolstoy,
{\it Rep. Math. Phys.} {\bf 4} 97 -- 111
\medskip

\noindent
[6].~H. de Guise and R. T. Sharp,
{\it J. Phys. A: Math. Gen.}
{\bf 24} 557 -- 568
\medskip

\noindent
[7].~{\v C}. Burdik, C. J. Cummins, R. W. Gaskell and
R. T. Sharp, {\it J. Phys. A: Math. Gen.}
4835 -- 4846
\medskip

\noindent
[8].~L. Farell, C. S. Lam and R. T. Sharp,
{\it J. Phys. A: Math. Gen.} {\bf 27} 2761 -- 2771
\medskip

\noindent
[9].~N. Hambli and R. T. Sharp,
{\it J. Phys. A: Math. Gen.} {\bf 28} 2581 -- 2588
\medskip

\noindent
[10].~R. Gaskell, A. Peccia and R. T. Sharp,
{\it J. Math. Phys.} {\bf 19} 727 -- 733
\medskip

\noindent
[11].~J. Patera and R. T. Sharp,
1979 Generating
Functions for Characters of Group Representations
and Their Applications Lecture Notes in Physics
{\bf 94} 175 -- 183
\medskip

\pagebreak

\end{document}